%&amstex          
\input amstex\documentstyle{amsppt}  
\pagewidth{12.5cm}\pageheight{19cm}\magnification\magstep1
\topmatter
\title Bruhat decomposition and applications\endtitle
\author G. Lusztig\endauthor
\address{Department of Mathematics, M.I.T., Cambridge, MA 02139}\endaddress
\thanks{Supported in part by the National Science Foundation}\endthanks
\dedicatory{Lecture given at a one day conference in memory of F. Bruhat held 
at the Institut Poincar\'e, Paris, on June 22, 2010.}\enddedicatory
\endtopmatter   
\document

\define\dw{\dot w}

\define\uuG{\un{\un G}}

\define\uW{\un W}

\define\sqc{\sqcup}

\define\part{\partial}
\define\em{\emptyset}

\define\m{\mapsto}

\define\sub{\subset}    

\define\T{\times}

\define\nl{\newline}

\define\un{\underline}

\define\g{\gamma}

\define\s{\sigma}

\define\Ph{\Phi}

\define\kk{\bold k}

\define\CC{\bold C}

\define\NN{\bold N}

\define\RR{\bold R}

\define\cb{\Cal B}

\define\AO{A}
\define\BRI{B1}
\define\BR{B2}
\define\BT{BT}
\define\BO{Bo}
\define\BTI{TB}
\define\DL{DL}
\define\HC{H}
\define\CHT{C1}
\define\CH{C2}
\define\EH{E}
\define\IW{I}
\define\IM{IM}
\define\GN{GN}
\define\MA{M}
\define\RO{R}
\define\SH{Sh}
\define\ST{S}
\define\LU{L}

\subhead 1. Statement\endsubhead
Let $\kk$ be a field. We assume that $\kk$ is algebraically closed, unless otherwise specified. Let G be a 
connected reductive algebraic group over $\kk$. Let $\cb$ be the variety of Borel subgroups of $G$. Let $B\in\cb$,
let $T$ be a maximal torus of $B$ and let $N$ be the normalizer of $T$ in $G$. Let $W=N/T$ be the Weyl group. For
any $w\in W$ let $G_w=B\dw B$ where $\dw\in N$ represents $w$. The following is a restatement of Theorem 7.1 in 
Bruhat's 1956 paper \cite{\BR}.

($*$) {\it The sets $G_w$ (with $w\in W$) form a partition of $G$.}
\nl
(Actually in \cite{\BR} it is assumed that $\kk=\CC$ and that $G$ is semisimple; the name "Borel subgroup" is not
used in \cite{\BR}. A variant of ($*$) over real numbers is also considered in \cite{\BR}.)

The partition $G=\sqc_{w\in W}G_w$ has been called the "Bruhat decomposition" in the Chevalley Seminar 
\cite{\CH, p.148}. In this talk we will examine the history (see no.2) and applications (see no.3,4) of this 
decomposition.

\subhead 2. History\endsubhead
In 1809, Gauss introduced his elimination method for solving systems of linear equations. As a consequence, 
almost any matrix in $G=GL_n(\CC)$ is a product $LU$ where $L$ (resp. $U$) is a lower (resp. upper) triangular 
matrix. The consideration of the subset $LU$ of $G$ is equivalent (up to translation) to the consideration of the
piece $G_{w_0}$ ($w_0$ being the element of $W$ of maximal length for the standard length function $l:W@>>>\NN$) 
in the Bruhat decomposition of $G$.

Actually, the Gauss elimination method and the $LU$ decomposition already appear in Chapter 8 of the chinese 
classic "The Nine Chapters on the Mathematical Art" (from the 2nd century, Han dynasty). Here the method of 
solving systems of linear equations is presented by means of examples involving systems with up to five unknowns.
(This is also the first place where negative numbers appear in the literature.)

In his 1934 paper \cite{\EH}, Ehresmann shows that any partial flag manifold of $G=GL_n(\CC)$ admits a 
decomposition into finitely many complex cells, generalizing earlier work of Schubert (around 1880) for the 
Grassmannians. The decomposition is in terms of a fixed full flag and the pieces of the decomposition are clearly
stable under the stabilizer of that fixed flag. Ehresmann's decomposition of the full flag manifold can now be 
viewed as induced by the Bruhat decomposition. Ehresmann also parametrizes his cells in terms of certain tableaux
of integers (for the full flag manifold of $GL_4(\CC)$ he describes explicitly the $4!$ tableaux which appear) 
but he does not interpret these tableaux in terms of the Weyl group.

In his 1951 paper \cite{\ST} (submitted in 1949), Steinberg identifies the orbits of $GL_n(F_q)$ on pairs of 
complete flags in $F_q^n$ with permutations of $n$ objects (see p.275,276), a result very close to $(*)$ for 
$GL_n$ over a finite field.

In their 1950 book \cite{\GN, p.122}, Gelfand and Naimark state and prove $(*)$ for $G=SL_n(\CC)$.

In his 1954 announcement \cite{\BRI}, Bruhat formulates for the first time ($*$) for general semisimple groups 
over $\CC$ and states that he has verified it for all classical groups. 

One of the consequences of ($*$) is that (when $\kk=\CC$), $\cb$ has no odd integral homology and no torsion in 
even integral homology. This was first proved by Bott \cite{\BO} in 1954 (independently of ($*$)) using Morse 
theory.

A proof of ($*$) (with $\kk=\CC$) valid for any $G$ was given in the 1956 paper \cite{\HC} of Harish-Chandra; 
this is the proof reproduced in Bruhat's 1956 paper \cite{\BR}. A proof of ($*$) for arbitrary $\kk$ was given in 
Chevalley's 1955 paper \cite{\CHT}; in this paper Chevalley mentions the existence of Harish-Chandra's proof 
(which was unpublished at the time). In \cite{\BTI}, Borel and Tits proved a version of ($*$) valid over an 
arbitrary field.

\subhead 3. Significance\endsubhead
By allowing one to reduce many questions about $G$ to questions about the Weyl group $W$, Bruhat 
decomposition is indispensable for the understanding of both the structure and representations of $G$. 
We shall illustrate this by several examples. (A further example is given in no.4.)

The order of a Chevalley group over a finite field was computed in \cite{\CHT} (using Bruhat decomposition) in 
terms of the exponents of the Weyl group. 

In the representation theory of a reductive group over a finite field a key role is played by the Iwahori-Hecke 
algebra (introduced in \cite{\IW}); this is a deformation of the group algebra of the Weyl group whose definition
is based on the Bruhat decomposition. In the same theory the varieties (introduced in \cite{\DL}) obtained by 
taking inverse image of a Bruhat double coset under the Lang map have turned out to be very useful.

In the representation theory of complex reductive groups, 
to understand the character of irreducible representations one needs the local intersection cohomology of 
the closure of a Bruhat double coset. 

In the representation theory of split $p$-adic reductive groups, a key role is played by the affine Hecke algebra
whose definition is based on the generalization of the Bruhat decomposition given by Iwahori and Matsumoto 
\cite{\IM}. (The non-split case was treated by Bruhat and Tits in \cite{\BT}.)

The folowing is a reformulation of ($*$):

(i) the orbits of $G$ acting on $\cb\T\cb$ by simultaneous conjugation are in natural bijection with $W$.
\nl
Assume now that $G$ is semisimple, simply connected and $\kk=\CC$. Note that (i) has been extended in two 
different directions as follows.

(ii) Let $G_\RR$ be the group of real points for a fixed real structure on $G$. In 1966 Aomoto \cite{\AO} showed 
that that the conjugation action of $G_\RR$ on $\cb$ has finitely many orbits and in 1979 Rossmann \cite{\RO} 
gave a parametrization of the orbits in terms of Weyl groups. 

(iii) Let $K$ be the identity component of the group of fixed points of an involution $\s$ of $G$ (as an 
algebraic group). In 1979 Matsuki \cite{\MA} showed that the conjugation action of $K$ on $\cb$ has finitely many
orbits and gave an explicit parametrization of the orbits (they are in bijection with a set of orbits as in (ii)).

The local intersection cohomology of the orbit closures in (iii) plays a key role for understanding the character
of irreducible representations of $G_\RR$ as in (ii).

\subhead 4. Bruhat decomposition and conjugacy classes\endsubhead
Assume that $G$ is semisimple and that the characteristic of $\kk$ is not a bad prime for $G$. By studying the 
interaction of Bruhat decomposition with conjugacy classes in $G$ we obtain a surprising connection between the 
set $\uuG$ of unipotent conjugacy classes in $G$ and the set $\uW$ of conjugacy classes in $W$. 

Let $C\in\uW$, let $d_C=\min(l(w);w\in C)$ and let $C_{min}=\{w\in C;l(w)=d_C\}$. We have the following result, 
see (\cite{\LU, 0.4(i)}):

(a) {\it Let $w\in C_{min}$. There is a unique $\g\in\uuG$ such that $\g\cap G_w\ne\em$ and such that whenever 
$\g'\in\uuG$, $\g'\cap G_w\ne\em$, we have $\g\sub\bar\g'$. Moreover, $\g$ depends only on $C$, not on $w$; we 
denote it by $\g_C$.}
\nl
Thus we obtain a map $\Ph:\uW@>>>\uuG$, $C\m\g_C$. 

Let $\uW_{el}$ be the set of conjugacy classes in $W$ which are elliptic (that is consist of elements with no 
eigenvalue $1$ in the reflection representation of $W$). Here are some properties of $\Ph$.

(b) $\Ph$ is surjective;

(c) $\Ph|_{\uW_{el}}:\uW_{el}@>>>\uuG$ is injective.

(d) If $C\in\uW_{el}$ and $w\in C_{min}$, then $\Ph(C)$ is the unique unipotent class $\g$ of $G$ such that 
$\g\cap G_w$ is a union of finitely many $B$-orbits for the conjugacy action 
of $B$ on $G_w$. Moreover, if $g\in\Ph(C)\cap G_w$, 
the dimension of the centralizer of $g$ in $B$ (resp. $G$) is equal to $0$ (resp. $d_C$).

(e) If $C\in\uW-\uW_{el}$, then $\Ph(C)$ has a simple description in terms of the map analogous to $\Ph$ for a 
Levi subgroup of a proper parabolic subgroup of $G$.
\nl
There is substantial evidence that the definition and properties of $\Ph$ are valid without assumption on the 
characteristic of $\kk$ and that the first assertion of (d) is valid for any $C\in\uW$.

For example, if $G=GL_n(\kk)$ then both $\uW$ and $\uuG$ may be identified with the set $P_n$ of partitions of 
$n$ (using the cycle types for $\uW$ and the sizes of the Jordan blocks for $\uuG$) and $\Ph$ becomes the 
identity map.

If $G=Sp_{2n}(\kk)$ then $W$ can be naturally identified with a subgroup of the symmetric group in $2n$ letters 
hence we have a natural map $i:\uW@>>>P_{2n}$ (neither injective nor surjective in general); also, via the 
obvious imbedding $G\sub GL_{2n}(\kk)$ we obtain an imbedding of $\uuG$ into the set of unipotent classes of 
$GL_{2n}(\kk)$ hence we have a natural injective map $j:\uuG@>>>P_{2n}$. It turns out that the image of $i$ 
coincides with the image of $j$ and $\Ph:\uW@>>>\uuG$ is characterized by $j\Ph(C)=i(C)$ for all $C$.


\widestnumber\key{\GN}  
\Refs   
\ref\key{\AO}\by K.Aomoto\paper On some double coset decompositions of complex semisimple groups\jour J. Math. 
Soc. Japan\vol18\yr1966\pages1-44\endref
\ref\key{\BTI}\by A.Borel and J.Tits\paper Groupes r\'eductifs\jour Publ. Math. IHES\vol27\yr1965\pages55-152
\endref
\ref\key{\BO}\by R.Bott\paper On torsion in Lie groups\jour Proc.Nat.Acad.Sci.\vol40\yr1954/\pages586-588\endref
\ref\key{\BRI}\by F.Bruhat\paper Representations induites des groupes de Lie semisimples complexes\jour Comptes 
Rendues Acad. Sci. Paris\vol238\yr1954\pages437-439\endref
\ref\key{\BR}\by F.Bruhat\paper Sur les representations induites des groupes de Lie\jour Bull. Soc. Math. France
\vol84\yr1956\pages97-205\endref
\ref\key{\BT}\by F.Bruhat and J.Tits\paper Groupes r\'eductifs sur un corps local,I\jour Publ. Math. IHES\vol41
\yr1972\pages5-276\endref
\ref\key{\CHT}\by C.Chevalley \paper Sur certains groups simples\jour Tohoku Math.J.\vol7\yr1955\endref
\ref\key{\CH}\by C.Chevalley \book Classification des groupes alg\'ebriques semisimples\publ Springer\endref
\ref\key{\DL}\by P.Deligne and G.Lusztig\paper Representations of reductive groups over finite fields\jour Ann.
math.\vol103\yr1976\pages103-161\endref
\ref\key{\EH}\by C.Ehresmann\paper Sur la topologie de certains espaces homogenes\jour Ann. Math.\yr1934\endref
\ref\key{\GN}\by I.M.Gelfand and M.A.Naimark\book Unitarnye predstavleniya klassiceskih grupp\bookinfo Trudy Mat.
Inst. Steklov\vol36\publaddr Moscow\yr1950\endref
\ref\key{\HC}\by Harish-Chandra\paper On a lemma of Bruhat\jour J.Math.Pures Appl.\vol35\yr1956\pages203-210
\endref
\ref\key{\IW}\by N.Iwahori\paper On the structure of a Hecke ring of Chevalley groups over a finite field\jour 
J. Fac. Sci. Univ.Tokyo\vol10\yr1964\endref
\ref\key{\IM}\by N.Iwahori and H.Matsumoto\paper On some Bruhat decomposition and the structure of the Hecke 
rings of $p$-adic Chevalley groups\jour Publ. Math. IHES\vol25\yr1965\pages5-48\endref
\ref\key{\LU}\by G.Lusztig\paper From conjugacy classes in the Weyl group to unipotent classes\jour
arxiv:1003.0412\endref
\ref\key{\MA}\by T.Matsuki\paper The orbits of affine symmetric spaces under the action of minimal parabolic 
subgroups\jour J. Math. Soc. Japan\vol31\yr1979\pages331-357\endref
\ref\key{\RO}\by W.Rossmann\paper The structure of semisimple symmetric spaces\jour Canad. J. Math.\vol31\yr1979
\pages157-180\endref
\ref\key{\SH}\by K.Shen\book The Nine Chapters on the Mathematical Art\publ Oxford Univ. Press\yr1999\endref
\ref\key{\ST}\by R.Steinberg\paper A geometric approach to the representations of the full linear group over a 
Galois field\jour Trans. Amer. Math. Soc.\vol71\yr1951\pages274-282\endref
\endRefs
\enddocument